\documentclass[12pt]{article}
\usepackage{color}

\usepackage[italian, english]{babel}

\usepackage{amssymb,amsmath} 

\usepackage[OT2, T1]{fontenc} 

\usepackage[utf8]{inputenc}

\allowdisplaybreaks

\usepackage{hyperref}
\hypersetup{hidelinks}

\usepackage{scalerel}
\usepackage{tikz}

\usetikzlibrary{svg.path}
\definecolor{orcidlogocol}{HTML}{A6CE39}
\tikzset{
 orcidlogo/.pic={
 \fill[orcidlogocol] svg{M256,128c0,70.7-57.3,128-128,128C57.3,256,0,198.7,0,128C0,57.3,57.3,0,128,0C198.7,0,256,57.3,256,128z};
 \fill[white] svg{M86.3,186.2H70.9V79.1h15.4v48.4V186.2z}
 svg{M108.9,79.1h41.6c39.6,0,57,28.3,57,53.6c0,27.5-21.5,53.6-56.8,53.6h-41.8V79.1z M124.3,172.4h24.5c34.9,0,42.9-26.5,42.9-39.7c0-21.5-13.7-39.7-43.7-39.7h-23.7V172.4z}
 svg{M88.7,56.8c0,5.5-4.5,10.1-10.1,10.1c-5.6,0-10.1-4.6-10.1-10.1c0-5.6,4.5-10.1,10.1-10.1C84.2,46.7,88.7,51.3,88.7,56.8z};
 }
}
\newcommand\orcidicon[1]{\href{https://orcid.org/#1}{\mbox{\scalerel*{
\begin{tikzpicture}[yscale=-1,transform shape]
\pic{orcidlogo};
\end{tikzpicture}
}{|}}}}

\newtheorem{theorem}{Theorem}[section]
\newtheorem{proposition}[theorem]{Proposition}

\newtheorem{remark}[theorem]{Remark}
\newtheorem{definition}[theorem]{Definition}

\def\claim#1.{\noindent {\bf #1.}}

\def\bull{\vrule height 1.8ex width 1.0ex depth .1ex }
\def\flushright#1{{\unskip\nobreak\hfil\penalty50\hskip2em\hbox{}\nobreak\hfil #1\parfillskip=0pt\finalhyphendemerits=0\par}}
\def\QED{\ifmmode\eqno\hbox{$\bull$}\else\flushright{\hbox{$\bull$}}\fi}

\newcommand{\norm}[1]{\left\Vert#1\right\Vert}
\newcommand{\tnorm}[1]{{\left\vert\kern-0.25ex\left\vert\kern-0.25ex\left\vert #1 
		\right\vert\kern-0.25ex\right\vert\kern-0.25ex\right\vert}}
		
\newcommand{\mc}[1]{\mathcal{#1}}

\newcommand{\parag}[1]{\left\{ \begin{aligned} #1 \end{aligned}\right.}

\def\cat{\mathop{\rm cat}\nolimits}
\newcommand{\cupl}{{\rm cupl}}

\newcommand{\R}{\mathbb{R}}
\newcommand{\N}{\mathbb{N}}
\newcommand{\F}{\mathbb{F}}

\newcommand{\eps}{\varepsilon}

\begin{document}
	

	\title{Multiplicity and concentration results for local and fractional NLS equations with critical growth}

	\author{
		Marco Gallo \orcidicon{0000-0002-3141-9598}
		\\Dipartimento di Matematica
		\\Universit\`{a} degli Studi di Bari Aldo Moro
		\\Via E. Orabona 4, 70125 Bari, Italy
		\\ Email: marco.gallo@uniba.it
		\\ 
	}

	\maketitle

\abstract{
Goal of this paper is to study positive semiclassical solutions of the nonlinear Schr\"odinger equation 
$$ \eps^{2s}(- \Delta)^s u+ V(x) u= f(u), \quad x \in \R^N,$$
where $s \in (0,1)$, $N \geq 2$, $V \in C(\R^N,\R)$ is a positive potential and $f$ is assumed critical and satisfying general Berestycki-Lions type conditions. 
We obtain existence and multiplicity for $\eps>0$ small, where the number of solutions is related to the cup-length of a set of local minima of $V$. Furthermore, these solutions are proved to concentrate in the potential well, exhibiting a polynomial decay. 
We highlight that these results are new also in the limiting local setting $s=1$ and $N\geq 3$, with an exponential decay of the solutions.
}

	\bigskip
	
\textbf{Keywords:}
Fractional Laplacian, Nonlinear Schr\"odinger Equation, Critical Exponent, Singular Perturbation, Spike Solutions, Cup-length

\medskip

\textbf{MSC:}
 35A15, 35B25, 35B33, 35Q55, 35R11, 47J30, 58E05

	\tableofcontents
	

\section{Introduction}\label{section:1}
Following a suggestion by Dirac, in 1948 Feynman proposed a new suggestive description of the time evolution of the state of a non-relativistic quantum particle. According to Feynman, the wave function solution of the Schr\"odinger equation should be given by a heuristic integral over the space of paths: the classical notion of a single, unique classical trajectory for a system is replaced by a functional integral over an infinity of quantum-mechanically possible trajectories. Following Feynman’s path integral approach to quantum mechanics, Laskin \cite{Las} developed a new extension of the fractality concept by replacing the path integral over Brownian motions (random motion seen in swirling gas molecules) with Lévy flights (a mix of long trajectories and short, random movements found in turbulent fluids), deriving the
fractional nonlinear Schr\"odinger (fNLS for short) equation 
\begin{equation}\label{fNLS}
i \hbar \partial_t \psi = \hbar^{2s} (-\Delta)^s \psi + V(x) \psi- f(\psi), \quad (t,x) \in (0,+\infty) \times \R^N.
\end{equation}
Here $s \in (0,1)$, $N > 2s$, the symbol $(-\Delta)^{s}\psi=\mc{F}^{-1}(|\xi|^{2s} \mc{F}(\psi))$ denotes the fractional power of the Laplace operator defined via Fourier transform $\mc{F}$ on the spatial variable, $\hbar$ designates the usual Planck constant, $V$ is a real potential and $f$ is a Gauge invariant nonlinearity, i.e. $f(e^{i \theta} \rho) = e^{i \theta} f(\rho)$ for any $\rho, \theta \in \R$. The wave function $\psi(x,t)$ represents the quantum mechanical probability amplitude for a given unit mass particle to have position $x$ and time $t$, under the confinement due to the potential $V$. We refer to \cite{Las, Las2, Las3} for a detailed discussion on the physical motivation of the {fNLS} equation, and we highlight that several applications in the physical sciences could be mentioned, ranging from the description of boson stars to water wave dynamics, from image reconstruction to jump processes in finance.

Special solutions of the equation \eqref{fNLS} are given by the standing waves, i.e. factorized functions $\psi (t,x) = e^{\frac{i E t}{\hbar}} u(x)$ with $E\in \R$. 
For small $\hbar >0$, these standing waves are usually called \emph{semiclassical} states since the transition from quantum physics to classical physics is somehow described letting $\hbar \to 0$: roughly speaking, when $s = 1$ the centers of mass $q_{\eps}(t)$ of the soliton solutions in \eqref{fNLS}, under suitable assumptions and suitable initial conditions, converge as $\hbar\to 0$ to the solution of the Newton equation's of motion
\begin{equation}\label{eq_newton}
\ddot{q}(t)=-\nabla V(q(t)), \quad t \in (0, +\infty);
\end{equation}
for $s\in (0,1)$ a suitable power-type modification of equation \eqref{eq_newton} is needed. Here, considering small $\hbar$ roughly means that the size of the support of the soliton in \eqref{fNLS} is considerably smaller than the size of the potential $V$: for details we refer to \mbox{\cite{BrJe,
FGJS,JFG,BGM1,BGM2},} and to \cite{SS} for the fractional case.

Without loss of generality, shifting $E$ to
$0$ and denoting $\hbar \equiv \eps$, the search for semiclassical states leads to investigate the following nonlocal equation
\begin{equation}\label{eq_principale}
\eps^{2s} (-\Delta)^s u + V(x) u = f(u), \quad x \in \R^N
\end{equation}
where $V$ is positive and $\eps > 0$ is small. 
Setting $v:= u(\eps \cdot)$, we observe that \eqref{eq_principale} can be rewritten as
\begin{equation}\label{eq_cambio_var}
(-\Delta)^s v + V(\eps x) v= f(v), \quad \textnormal{ $x \in \R^N$},
\end{equation}
thus the equation
\begin{equation}\label{eq_limite} 
(-\Delta)^s U + m_0 U = f(U), \quad x \in \R^N
\end{equation}
becomes a formal limiting equation, as $\eps \to 0$, of \eqref{eq_cambio_var}. 
Solutions of \eqref{eq_principale} usually exhibit concentration phenomena as $\eps\to 0$: by \emph{concentrating solutions} we mean a family $u_{\eps}$ of solutions of \eqref{eq_principale} which converges, up to rescaling, to a ground state of \eqref{eq_limite} and whose maximum points converge to some point $x_0 \in \R^N$ given by the topology of $V$ (see Theorem \ref{teo_main} for a precise statement). This point $x_0$ reveals, most of the time, to be a critical point of $V$ - i.e. an equilibrium of \eqref{eq_newton} - as shown by \cite{Wan, FMV}.

In the subcritical case (see (f4) below), the semiclassical analysis of local NLS equations has been largely investigated 
starting from the seminal papers \cite{FW, Oh2}: here the authors implement a Lyapunov-Schmidt dimensional reduction argument to gain existence of solutions for homogeneous sources, relying on the nondegeneracy of the ground states of the limiting problem \eqref{eq_limite}. Successively, variational techniques have been implemented to gain both existence and multiplicity, see \cite{Rab, Wan, ABC, DF, CL, AMS, BJ3, BT1, CJT} and references therein.

In the last years the subcritical fractional case has gradually aroused interest in many authors, and we confine to mention \cite{DPW, FMV, AlMi1, FS, Seo, Che}; in particular we refer to the recent paper \cite{CG} and the references therein. 
In \cite{CG} the author and Cingolani gain multiplicity and concentration of solutions assuming the following set of assumptions: on $V$ they assume
\begin{itemize}
\item[(V1)] $\, V\in C(\R^N, \R)\cap L^{\infty}(\R^N)$, $\underline{V}:=\inf_{\R^N} V>0$ (see also Remark \ref{rem_ipotesi_V}),
\item[(V2)] $\,$there exists a bounded domain $\Omega\subset \R^N$ such that
$$m_0:= \inf_{\Omega} V < \inf_{\partial \Omega} V,$$
with set of local minima
\begin{equation}\label{eq_def_K}
K:=\{ x \in \Omega \mid V(x)=m_0\},
\end{equation}
\end{itemize}
while on $f$ they assume Berestycki-Lions \cite{BL} type conditions, i.e.
\begin{itemize}
\item[(f1)] $f\in C^{0,\gamma}_{loc}(\R, \R)$ for some $\gamma \in (1-2s, 1)$ if $s\in (0,1/2]$ (see also Remark \ref{rem_ipotesi_reg}),
\item[(f2)] $f(t)\equiv 0$ for $t \leq 0$,
\item[(f3)] $\, \lim_{t \to 0} \frac{f(t)}{t}=0$,
\item[(f4)] $\, \lim_{t \to +\infty} \frac{f(t)}{|t|^q}=0$ for some $q \in (1, 2^*_s-1)$, where $2^*_s:=\frac{2N}{N-2s}$ is the Sobolev critical exponent,
\item[(f5)] $\, F(t_0)> \frac{1}{2} m_0 t_0^2$ for some $t_0>0$, where $F(t):= \int_0^t f(\tau) d\tau$.
\end{itemize}

For $s \in (0,1)$ the existence, in the critical setting, of a solution under (V1)-(V2) has been faced by \cite{JLZ} with $V\in C^1(\R^N)$, and in \cite{He2}. Here it results crucial the use of $L^{\infty}$-estimates, which are nontrivial in the critical framework. 
The authors assume almost optimal Berestycki-Lions type assumptions on $f$, that is $f\in C^1(\R)$, (f2)-(f3) and
\begin{itemize}
\item[(f4')] $\lim_{t \to +\infty} \frac{f(t)}{t^{2^*_s-1}}=a>0$,
\item[(f5')] for some $C>0$ and $\max\{2^*_s-2, 2\} <p< 2^*_s$, i.e. satisfying
\begin{equation}\label{eq_cond_p}
p \in \parag{&\Big(\frac{4s}{N-2s}, \frac{2N}{N-2s}\Big)& \quad N \in (2s, 4s), \\ &\Big(2, \frac{2N}{N-2s}\Big)& \quad N\geq 4s, }
\end{equation}
(see also Remark \ref{rem_C_grande}), it results that
$$f(t) \geq t^{2^*_s-1} + C t^{p-1} \quad \textnormal{for $t \geq 0$}.$$
\end{itemize}
Notice that the stronger condition on $p$ in the first line of \eqref{eq_cond_p} is verified, whenever $N\geq 2$, only if $N=2$ and $s \in (\frac{1}{2}, 1]$, or $N=3$ and $s\in (\frac{3}{4}, 1]$. 
We point out that the condition $C>0$ in (f4') is of key importance: indeed, for pure critical nonlinearities of the type $f(t)=|t|^{2^*_s-t}t$, the limiting problem \eqref{eq_limite} exhibits a quite different scenario. 

Inspired by \cite{Rab}, multiplicity of solutions in the case of global minima of $V$ was studied in \cite{SZY} for power-type nonlinearities, and 
in \cite{LTZW}: in this paper the authors gain existence and multiplicity for functions of the type
\begin{equation}\label{eq_cond_g}
f(t)=g(t)+ |t|^{2^*_s-2}t,
\end{equation}
where $g$ is subcritical and it satisfies a monotonicity condition. Thanks to this last assumption, the Nehari manifold tool can be implemented, and the number of solutions is related to the category of the set of global minima.

Existence of multiple solutions for local minima has been investigated, in the spirit of \cite{DF}, by \cite{HeZo
}, with sources of the type \eqref{eq_cond_g} where now $g$ satisfies also an Ambrosetti-Rabinowitz condition: this assumption enables to employ mountain pass and Palais-Smale arguments, combined with a penalization scheme. Again, the authors are able to find $\cat(K)$ solutions, where $K$ is the set of local minima of $V$.

\smallskip

In the present paper we extend the results in \cite{CG} to the critical case, proving multiplicity of positive solutions for the fNLS equation \eqref{eq_principale} when $\eps$ is small, without assuming monotonocity and Ambrosetti-Rabinowitz conditions on $f$, nor nondegeneracy and global conditions on $V$.

We prove the following theorem. Up to the author's knowledge, this result is new and improve the results in \cite{
HeZo, He2, LTZW}.
\begin{theorem}\label{teo_main}
Suppose $s \in (0,1)$, $N\geq 2$ and that \textnormal{(V1)-(V2)}, \textnormal{(f1)--(f3)}, \textnormal{(f4')-(f5')} hold. 		
Let $K$ be defined by \eqref{eq_def_K}.
Then, up to a discretized subsequence, for small $\eps>0$ equation \eqref{eq_principale} has at least $\cupl(K)+1$ positive solutions, which belong to $C^{0, \sigma}(\R^N) \cap L^{\infty}(\R^N)$ for some $\sigma \in (0,1)$. Moreover, each of these sequences $u_{\eps}$ concentrate in $K$ as $\eps \to 0$. Namely, for each small $\eps >0$ there exists a maximum point $x_{\eps}\in \R^N$ of $u_{\eps}$ such that 
$$\lim_{\eps \to 0}d(x_{\eps}, K) =0.$$
In addition, $ u_{\eps}(\eps \cdot+x_{\eps})$ converges in $H^s(\R^N)$ and uniformly on compact sets to a least energy solution of \eqref{eq_limite} and, for some positive $C', C''$ independent on $\eps$, we have the uniform polynomial decay
$$\frac{C'}{1+|\frac{x-x_{\eps}}{\eps}|^{N+2s}}\leq u_{\eps}(x) \leq \frac{C''}{1+|\frac{x-x_{\eps}}{\eps}|^{N+2s}} \quad \textit{ for $x \in \R^N$}.$$
\end{theorem}

Here $\cupl(K)$ denotes the \emph{cup-length} of $K$ defined by the Alexander-Spanier cohomology with coefficients in some field $\F$ (see Definition \ref{def_cuplength}). This topological tool denotes the geometric complexity of the set $K$, and it was successfully implemented also in \cite{AMS, CJT, Che, CG}.

\begin{remark}\label{examples}
Notice that the cup-length of a set $K$ is strictly related to the \emph{category} of $K$. Indeed, if $K$ is a contractible set (e.g. a point) or it is finite, then $\cupl(K)+1=\cat(K)=1$; if $K=S^{N-1}$ is the $N-1$ dimensional sphere in $\R^N$, then $\cupl(K) + 1= \cat(K) =2$; if $K=T^N$ is the $N$-dimensional torus, then $\cupl(K) + 1 = \cat(K)= N + 1$.
However in general $\cupl(K) +1 \leq \cat(K)$.
\end{remark}

The idea of the paper is the following: after having gained compactness and uniform $L^{\infty}$-bounds on the set of ground states of \eqref{eq_principale}, we introduce a suitable truncation on the nonlinearity $f$. The new truncated function
reveals to be subcritical, i.e. satisfying (f1)--(f5). Therefore, in the spirit of \cite{CG}, we employ a penalization argument on a neighborhood of expected solutions, perturbation of the ground states of a limiting problem, and this neighborhood reveals to be invariant under the action of a deformation flow. Compactness is restored also by the use of a new fractional center of mass, which engage a seminorm stronger than the usual Gagliardo one; the topological machinery between two level sets of the associated indefinite energy functional is then built also through the use of a Pohozaev functional. The number of solutions is thus related to the cup-length of $K$ and these solutions are proved to exhibit a polynomial decay and to converge to a ground state of the limiting equation. 
This last convergence allows finally to prove that these solutions solve the original problem \eqref{eq_limite}.

\begin{remark}\label{rem_ipotesi_V}\label{rem_C_grande} \label{rem_ipotesi_reg}
As observed in \cite{CG, CJT}, assumption \textnormal{(V1)} in Theorem \ref{teo_main} can be relaxed without assuming the boundedness of $V$ (see also \cite{BJ3, BT1}). 
Moreover, the condition $p>\max\{2^*_s-2, 2\}$ in \textnormal{(f5')} can be relaxed in $p>2$ by paying the cost of considering a sufficiently large $C\gg 0$; see for instance \cite{SZ2, He2}. 
Finally, we remark that \textnormal{(f1)}, instead of the mere continuity of $f$, is needed only to get a Pohozaev identity by means of the regularity of solutions (see \cite[Proposition 1.1]{BKS}). 
\end{remark}

We highlight that the conclusions of Theorem \ref{teo_main} hold also for $s=1$ and $N\geq 3$, see Theorem \ref{teo_main_locale}. In this local framework, previous results were given by \cite{ADS, AlSo, ZhZo2, WLXF, Amb4}: in particular, we extend here the existence result in \cite{ZCZ} to a multiplicity result. 
 In this setting, the solutions decay exponentially and enjoy more regularity. 
 Notice that in such a case (f1) means $f$ merely continuous.

\medskip

The paper is organized as follows. In Section \ref{sec_recalls} we recall some notions on the fractional Sobolev space and on the fractional Laplacian, together with some crucial $L^{\infty}$-bound on the critical limiting problem. Then in Section \ref{sec_truncated} we introduce a truncation to bring the problem back to the subcritical case, and in Section \ref{sec_proof} we prove Theorem \ref{teo_main}. Finally in Section \ref{sec_locale} we deal with the local case.

\section{Some recalls}\label{sec_recalls}

Let $N\geq 2$ and $s \in (0,1)$. We define the Gagliardo seminorm
$$[u]_{s}^2:= \int_{\R^N} \int_{\R^N} \frac{|u(x)-u(y)|^2}{|x-y|^{N+2s}} \, dx \, dy$$
and the fractional Sobolev space \cite[Section 2]{DGV}
$$H^s(\R^N):=\left \{ u \in L^2(\R^N) \mid [u]_{s} < + \infty\right\}$$
endowed with the norm
$$\norm{u}_{H^s(\R^N)}^2 := \norm{u}_{2}^2 + [u]_{s}^2, \quad u \in H^s(\R^N);$$
here $\norm{\cdot}_{q}$ denotes the $L^q(\R^N)$-Lebesgue norm for $q \in [1, +\infty]$.
Moreover, for every $u \in H^s(\R^N)$ we define the fractional Laplacian
\begin{equation}\label{eq_def_frac_lapl}
(-\Delta)^s u(x):= C(n,s) \int_{\R^N} \frac{u(x)-u(y)}{|x-y|^{N+2s}} \, dx
\end{equation}
where the integral is in the principal value sense and $C(n,s)>0$. 

\medskip

Let us recall some crucial results on the limiting critical problem \eqref{eq_limite}. Define the energy $C^1$-functional $\mc{L}: H^s(\R^N) \to \R$
$$\mc{L}(U):=\frac{1}{2} \norm{U}_{H^s(\R^N)}^2 - \int_{\R^N} F(U) \, dx, \quad U \in H^s(\R^N)$$
related to the equation \eqref{eq_limite}, and the set
$$\widehat{S}:= \big\{ U \in H^s(\R^N)\setminus \{0\} \mid \textnormal{$U$ ground state solution of \eqref{eq_limite}, $ U(0)=\max_{\R^N} U$}\big\},$$
where by ground state solution we mean that each $U\in \widehat{S}$ solves $\mc{L}'(U)=0$ and its energy attains the least energy of $\mc{L}$, i.e.
$$\mc{L}(U)= E_{m_0},$$
where
$$E_{m_0}:=\inf \big\{ \mc{L}(V) \mid V\in H^s(\R^N)\setminus\{0\}, \; \mc{L}'(V)=0\big\}.$$
We state now some results on the set of ground states.

\begin{proposition}
Every $U\in \widehat{S}$ is positive and satisfies the Pohozaev identity, i.e.
\begin{equation}\label{eq_pohozaev}
\frac{N-2s}{2} \int_{\R^N} |(-\Delta)^{s/2}U|^2 \, dx - N \int_{\R^N} \left(F(U)-\frac{m_0}{2} U^2 \right)\, dx=0.
\end{equation}
\end{proposition}

\claim Proof.
The positivity is a straightforward consequence of assumption (f2). Moreover, the Pohozaev identity is gained by means of regularity results (see also Proposition \ref{prop_upperbound} below) and explicit computations on the $s$-harmonic extension problem; the arguments can be easily adapted from \cite[Proposition 1.1]{BKS} to the critical case. 
\QED

\smallskip

\begin{proposition}
The set $\widehat{S}$ is nonempty and compact in $H^s(\R^N)$. 
Moreover, every $U\in \widehat{S}$ is also a Mountain Pass solution for the problem \eqref{eq_limite}, i.e. set
$$C_{mp}:=\inf_{\gamma \in \Gamma} \sup_{t \in [0,1]} \mc{L}(\gamma(t))$$
with
$$\Gamma:=\big \{ \gamma \in C([0,1], \, H^s(\R^N)) \mid \gamma(0)=0, \; \mc{L}(\gamma(1))<0\big \}$$
we have
$$ E_{m_0}=C_{mp}.$$
Finally, up to a rescaling, every $U\in \widehat{S}$ is also a solution of the following constrained minimization problem
$$C_{min}:=\inf \big\{ \mc{T}(U) \mid \mc{V}(U)=1\big\}$$
where
$$\mc{T}(U):= \int_{\R^N} |(-\Delta)^{s/2}U|^2 \, dx, \quad \mc{V}(U):=\int_{\R^N} \left(F(U)-\frac{m_0}{2}U^2\right)\, dx;$$
in addition
\begin{equation}\label{eq_rel_ener_cmin}
E_{m_0}= \frac{s}{N} \left(\frac{N-2s}{2N}\right)^{\frac{N-2s}{2N}} (C_{min})^{\frac{N}{2s}}.
\end{equation}
\end{proposition}

\claim Proof.
Existence of a ground state solution can be achieved through the use of best Sobolev constants and minimization of $C_{min}$ as classically made by \cite{BL}; similarly, the equivalence with the mountain pass formulation is discussed as in \cite{JT1}. Arguing as in \cite[Lemma 3.4 and Proposition A.4]{CG} we instead gain compactness. We refer to \cite[Proposition 2.4
]{JLZ} for the precise statement and to \cite[Section 4.1 and Remark 1.2
]{ZMS} for details. 

Moreover, as observed in Remark \ref{rem_C_grande}, in order to show that $\widehat{S}\neq \emptyset$, the restriction on the range of $p$ in assumption (f5') can be substituted, by arguing as in \cite[Lemma 3.3]{SZY}, with the request that $C$ is sufficiently large (see also \cite{He2} and references therein).
\QED

\medskip

As a key property to employ the truncation argument we have the following result.
\begin{proposition}\label{prop_upperbound}
We have
$$\sup_{U \in \widehat{S}} \norm{U}_{\infty} < \infty.$$
\end{proposition}

\claim Proof.
An $L^{\infty}$-estimate for functions in $\widehat{S}$ is achieved by a Moser iteration argument \cite[Appendix B]{Str}; moreover, this estimate turns to be uniform thanks to the compactness of $\widehat{S}$, as done in \cite[Proposition A.4]{CG}. We refer to \cite[Proposition 3.1 and Remark 1.3]{JLZ} for the precise statement and details (see also \cite{DMV}).
\QED

\section{The truncated problem}\label{sec_truncated}

In virtue of Proposition \ref{prop_upperbound}, let
$$M:= \sup_{U \in \widehat{S}} \norm{U}_{\infty} +1$$
and set
$$k:= \sup_{t \in [0, M]} f(t) \in [0, +\infty);$$
moreover define the \emph{truncated} nonlinearity
$$f_k(t):= \min\{ f(t), k\}, \quad t \in \R.$$
We have the following properties on $f_k: \R \to \R$:
\begin{itemize}
\item $f_k(t) \leq f(t)$ for each $t \in\R$,
\item $f_k(t)=f(t)$ whenever $|t|\leq M$,
\item $f_k(U)=f(U)$ for every $U \in \widehat{S}$.
\end{itemize}	
Notice that the same relations hold also for $F$ and $F_k(t):= \int_0^t f_k(\tau)d \tau$.

We have that $f_k$ is subcritical, i.e. $f_k$ satisfies assumptions (f1)--(f5) with $q \in (1, 2^*_s-1)$ however fixed, as long as we can choose $t_0$ in (f5) such that $t_0 \in [0,M]$. But this is easily achieved: indeed, fixed a whatever $U\in \widehat{S}$, by the Pohozaev identity \eqref{eq_pohozaev} we have (notice that $(-\Delta)^{s/2}U$ cannot identically vanish)
$$N \int_{\R^N} \left(F(U)-\frac{m_0}{2} U^2 \right)\, dx=\frac{N-2s}{2} \norm{(-\Delta)^{s/2}U}_2^2 >0$$
and thus there exists an $x_0 \in \R^N$ such that
$$F(U(x_0))> \frac{m_0}{2} U(x_0)^2;$$
setting $t_0:=U(x_0)\in [0,M]$ we have the claim.

Consider now the truncated problem
\begin{equation}\label{eq_principale_k}
\eps^{2s} (-\Delta)^s u + V(x) u = f_k(u), \quad x \in \R^N
\end{equation}
and the corresponding limiting truncated problem
\begin{equation}\label{eq_limite_k}
(-\Delta)^s U + m_0 U = f_k(U), \quad x \in \R^N.
\end{equation}
Defined
$$\widehat{S}_k:= \big\{U \in H^s(\R^N)\setminus \{0\} \mid \textnormal{$U$ ground state solution of \eqref{eq_limite_k}, $ U(0)=\max_{\R^N} U$}\big\},$$
we have that the following key relation hold.
\begin{proposition}\label{prop_uguagl_insiem}
It results that 
$\widehat{S}=\widehat{S}_k$. Moreover, the least energy levels 
coincide.
\end{proposition}

\claim Proof.
Let us denote by $\mc{L}_k$, $\Gamma_k$, $\mc{V}_k$, $E_{m_0}^k=C_{mp}^k$, $C_{min}^k$ the quantities of the truncation problem analogous to the abovementioned ones of the critical problem.

First observe that, by $\mc{L}_k \geq \mc{L}$, we have $\Gamma_k \subset \Gamma$ and
\begin{equation}\label{eq_dis_mp}
C_{mp}^k \geq C_{mp};
\end{equation}
moreover for any $V \in \widehat{S}$ we have also $\mc{L}'_k(V)=0$, and hence
\begin{equation}\label{eq_dis_min}
\min_{V \in \widehat{S}} \mc{L}_k(V) \geq \min_{\mc{L}'_k(V)=0} \mc{L}_k(V) = E^k_{m_0}.
\end{equation}
Let now $U \in \widehat{S}$. We have by \eqref{eq_dis_mp} and \eqref{eq_dis_min}
$$C_{mp}^k \geq C_{mp} = \mc{L}(U) = E_{m_0}= \min_{V \in \widehat{S}} \mc{L}(V)= \min_{V \in \widehat{S}} \mc{L}_k(V) \geq E^k_{m_0}.$$
Therefore
$$\mc{L}_k(U)=\mc{L}(U) = C_{mp}^k = E_{m_0}^k$$
which, together with $\mc{L}'_k(U)=\mc{L}'(U)=0$, gives that $U \in \widehat{S}_k$. Hence $\widehat{S}\subset \widehat{S}_k$. 
As a further consequence we gain
\begin{equation}\label{eq_energ_ugual}
 E_{m_0}^k = E_{m_0}.
 \end{equation}
We show now that $\widehat{S}_k \subset \widehat{S}$. 
By \eqref{eq_energ_ugual}, \eqref{eq_rel_ener_cmin} and the analogous relation on the subcritical problem, we have
$$C_{min}^k=C_{min},$$
thus, 
by rescaling, it is sufficient to prove that every minimizer of $C_{min}^k$ is also a minimizer of $C_{min}$. Let thus $U$ be a minimizer for $C_{min}^k$, i.e. $\mc{T}(U)=C_{min}^k$ and $\mc{V}_k(U)=1$. Since $\mc{T}(U)=C_{min}$, it suffices to prove that $\mc{V}(U)=1$. By definition, we have 
$$\mc{V}(U)\geq \mc{V}_k(U)=1.$$
On the other hand, set $\theta:= (\mc{V}(U))^{\frac{1}{N}}$ we obtain, by scaling, that $\mc{V}(U(\theta \cdot))=1$ and thus
$$\mc{T}(U)=C_{min} \leq \mc{T}(U(\theta \cdot))= \theta^{-\frac{N+2s}{N}} \mc{T}(U)$$
from which we achieve
$$\mc{V}(U)\leq 1.$$
This concludes the proof.
\QED

\section{Proof of Theorem \ref{teo_main}}\label{sec_proof}

Before proving the main result, 
we first recall the definition of \emph{cup-length} (see e.g. \cite{Cha, FLRW, CLOT} and references therein).

\begin{definition}\label{def_cuplength}
Let $A$ be a topological space, and let $\mathbb{F}$ be a whatever field. 
Denote by 
$$H^*(A)=\bigoplus_{q\geq 0} H^q(A)$$
the Alexander-Spanier cohomology with coefficients in $\mathbb{F}$ (see \cite{Hat} and references therein). Let 
$$\smile: H^*(A)\times H^*(A)\to H^*(A)$$
be the cup-product. We define the \emph{cup-length} of $A$ as 
\begin{eqnarray*}
\cupl(A):= \max \big\{ l \in \N \mid \; &\exists \alpha_0 \in H^*(A), \; \exists \alpha_i \in H^{q_i}(A), \, q_i \geq 1, \textit{ for $i=1 \dots l$,}\\
& s.t. \; \alpha_0 \smile \alpha_1 \smile \dots \smile \alpha_l \neq 0 \, \textit{ in $H^*(A)$} \big\};&
\end{eqnarray*}
if such $l\in \N$ does not exist but $H^*(A)$ is nontrivial, we have $\cupl(A):=0$, otherwise we define $\cupl(A):=-1$.
\end{definition}

In the case $A$ is not connected, a slightly different definition (which makes the cup-length additive) can be found in \cite{Bar}. See Example \ref{examples} for some computation and comparison with the notion of the Lusternik-Schnirelmann category.

\medskip

\claim Proof of Theorem \ref{teo_main}.
We first look at the truncated problem \eqref{eq_principale_k}. Indeed, by \cite[Theorem 1.1 and Theorem 1.4]{CG} we obtain the existence of $\cupl(K)+1$ sequences of solutions of \eqref{eq_principale_k} satisfying the properties of Theorem \ref{teo_main} for $\eps>0$ small. 
We give here only an outline of the proof; to avoid cumbersome notation, we omit the dependence on the value $k$. 

Through a compact slight perturbation of the set $\widehat{S}_k$ (see \cite[Section 3.2]{CG}), still called $\widehat{S}_k$, we first define, for each $r>0$, a non-compact neighborhood of $\widehat{S}_k$
$$S(r):=\Big\{ u=U(\cdot-y)+\varphi \in H^s(\R^N) \mid U \in \widehat{S}_k, \; y \in \R^N, \; \varphi \in H^s(\R^N), \; \norm{\varphi}_{H^s(\R^N)}< r \Big\}.$$
To detect information on its elements we define a \emph{minimal radius} map $\widehat{\rho}: H^s(\R^N) \to \R_+$
\begin{eqnarray*}
&\widehat{\rho}(u):= \inf \left \{ \norm{u-U(\cdot-y)}_{H^s(\R^N)} \mid U \in \widehat{S}_k, \; y \in \R^N \right\}, \quad u \in H^s(\R^N) ,&\\
&u \in S(r) \implies \widehat{\rho}(u) < r,&
\end{eqnarray*}
and, for some suitable $\rho_0, R_0>0$, a \emph{barycentric} map $\Upsilon: S(\rho_0) \to \R^N$
\begin{eqnarray*}
&\Upsilon(u):=\frac{\int_{\R^N} y \, d(y,u) dy}{\int_{\R^N} d(y,u) dy}, \quad u \in S(\rho_0) ,&\\
&u=U(\cdot-y)+\varphi \in S(\rho_0) \implies |\Upsilon(u)-y|\leq 2R_0;&
\end{eqnarray*}
the \emph{density} map $d(y, u)$ appearing in the center of mass is defined by
$$d(y,u):= \psi \left(\inf_{U\in \widehat{S}_k} \norm{u-U(\cdot-y)}_{B_{R_0}(y)}\right), \quad (y, u) \in \R^N \times S(\rho_0),$$
where $\psi$ is a suitable cut-off function (see \cite[Lemma 3.7]{CG}) and the norm involved is a modification of the usual $H^s(B_{R_0}(y))$-norm, made through the use of a stronger seminorm which takes into account the tails of the functions, i.e.
$$\norm{u}_A^2:= \int_A u^2 \, dx + \int_{A} \int_{\R^N} \frac{|u(x)-u(y)|^2}{|x-y|^{N+2s}} \, dx \, dy, \quad u \in H^s(\R^N), \; A \subset \R^N.$$
Then, in order to localize solutions in $\Omega$, we introduce a suitable penalization on the functional $I_{\eps}: H^s(\R^N) \to \R$
$$I_{\eps}(u):=\frac{1}{2} \norm{(-\Delta)^{s/2} u}_2^2 + \frac{1}{2} \int_{\R^N} V(\eps x) u^2 \, dx - \int_{\R^N}F_k(u) \, dx, \quad u \in H^s(\R^N) $$
associated with the rescaled equation
\begin{equation}\label{eq_cambio_var_k}
(-\Delta)^s u + V(\eps x) u= f_k(u), \quad \textnormal{ $x \in \R^N$},
\end{equation}
and we call this penalized functional $J_{\eps}: H^s(\R^N) \to \R$ (see \cite[Section 4.1]{CG}).
Then, we restrict our attention to a \emph{neighborhood of expected solutions}
$$\mathcal{X}_{\eps, \delta}:= \big\{ u \in S(\rho_0) \mid \eps \Upsilon(u) \in K_d, \; J_{\eps}(u) < E_{m_0} + R(\delta, \widehat{\rho}(u)) \big\},$$
where $K_d$ is a suitable neighborhood of $K$, $R(\delta, \widehat{\rho}(u))$ is a suitable $u$-dependent radius and $\delta>0$ is chosen sufficiently small (see \cite[Section 4.3]{CG}).
On $\mc{X}_{\eps, \delta}$, for $\eps$ small, we succeed in proving delicate $\eps$-independent gradient estimates for $J_{\eps}$, a truncated Palais-Smale-type condition, and the existence of a deformation flow $\eta: [0,1]\times \mc{X}_{\eps, \delta} \to \mc{X}_{\eps,\delta}$; moreover, we prove that each solution of $J_{\eps}'(u)=0$ is also a solution of the original problem $I_{\eps}'(u)=0$ (see \cite[Theorem 4.7, Corollary 4.9, Proposition 4.10 and Lemma 4.11]{CG}).

To find multiple solutions we build two continuous maps satisfying
$$I \times K \; \stackrel{\Phi_{\eps}} \to \; \mathcal{X}_{\eps, \delta}^{E_{m_0}+\hat{\delta}} \;\stackrel{\Psi_{\eps}} \to \; I \times K_d \quad \textnormal{ and }\quad \partial I \times K \;\stackrel{\Phi_{\eps}} \to \;\mathcal{X}_{\eps, \delta}^{E_{m_0}-\hat{\delta}} \;\stackrel{\Psi_{\eps}} \to\; (I\setminus \{1\}) \times K_d,$$
where $I\subset\R$ is a suitable neighborhood of $1$, $\widehat{\delta}\in(0,\delta)$ and the superscript denotes the intersection with the sublevels of $J_{\eps}$. These maps are defined by
\begin{eqnarray*}
&\Phi_{\eps}(t,y):= U_0\left(\tfrac{\cdot-y/\eps}{t}\right), \quad (t,y) \in I\times K,&\\
&\Psi_{\eps}(u):= \big( T(P_{m_0}(u)), \eps \Upsilon(u)\big), \quad u \in \mathcal{X}_{\eps, \delta}^{E_{m_0}+\hat{\delta}},&
\end{eqnarray*}
where $U_0\in \widehat{S}_k$ is fixed, $T$ is a truncation over the interval $I$, and $P_{m_0}$ is a Pohozaev functional such that $P_{m_0}(U)=1$ for every $U\in \widehat{S}_k$ (see \cite[Section 4.4]{CG}).
The composition $\Psi_{\eps} \circ \Phi_{\eps}$ results being homotopic to the identity, and this leads to the existence of at least $\cupl(K)+1$ solutions by the following chain of inequalities involving the relative category (see \cite[Section 5]{CG})
\begin{eqnarray*}
\lefteqn{\# \big\{ u \textnormal{ solutions of \eqref{eq_cambio_var_k}}\big\} 
\geq 
\# \Big\{ u \in (\mathcal{X}_{\eps, \delta})^{E_{m_0}+\hat{\delta}}_{E_{m_0}-\hat{\delta}} \mid J_{\eps}'(u)=0 \Big\}}
\\
 &\quad \geq& 
\cat\Big (\mathcal{X}_{\eps, \delta}^{E_{m_0}+\hat{\delta}},\, \mathcal{X}_{\eps, \delta}^{E_{m_0}-\hat{\delta}}\Big)
 \geq
 \cupl (\Psi_{\eps} \circ \Phi_{\eps}) +1\geq \cupl(K) +1.
\end{eqnarray*}
Finally, uniform $L^{\infty}$-bounds, $C^{0,\sigma}$-regularity and concentration in $K$ are proved through the use of recent fractional De Giorgi classes \cite{Coz} (see \cite[Section 5.1]{CG}).
In particular, rescaling back again, for each of these sequences $u_{\eps}$ of solutions of \eqref{eq_principale_k}, it is proved that there exist $C, C'>0$, $U$ ground state of \eqref{eq_limite_k} and $x_0 \in \R^N$ such that, called $x_{\eps}\in \R^N$ a maximum point of $u_{\eps}$, it results that (up to a discretized subsequence) 
$$\lim_{\eps\to 0}d(x_{\eps}, K) =0, \quad \frac{C'}{1+|\frac{x-x_{\eps}}{\eps}|^{N+2s}}\leq u_{\eps}(x) \leq \frac{C''}{1+|\frac{x-x_{\eps}}{\eps}|^{N+2s}} \quad \textit{ for $x \in \R^N$},$$
and
\begin{equation}\label{eq_conv_U}
u_{\eps}(\eps \cdot + x_{\eps}) \to U(\cdot + x_0) \in \widehat{S}_k
\end{equation}
in $H^s(\R^N)$ and locally on compact sets. For further details we refer to \cite{CG}. 

By Proposition \ref{prop_uguagl_insiem} we have $U(\cdot + x_0) \in \widehat{S}$. In particular, by the pointwise convergence in \eqref{eq_conv_U} we obtain
$$\norm{u_{\eps}}_{\infty}= u(x_{\eps})\to U(x_0) \leq \norm{U}_{\infty} < M$$
which implies
$$\norm{u_{\eps}}_{\infty} < M$$
definitely for $\eps$ small. As a consequence 
$$f_k(u_{\eps})=f(u_{\eps})$$
and thus $u_{\eps}$ are solutions of the original problem \eqref{eq_principale}, satisfying the desired properties.
\QED

\section{The local case}\label{sec_locale}

The arguments presented in Theorem \ref{teo_main} apply, with suitable modifications, also to local nonlinear Schr\"odinger equations. We give here some details. 
Conditions (f1)--(f3), (f4')-(f5') rewrite in the local case $s=1$ as
\begin{itemize}
\item[(f1'')] $f\in C(\R, \R)$,
\item[(f2'')] $f(t)\equiv 0$ for $t\leq 0$,
\item[(f3'')] $\, \lim_{t \to 0} \frac{f(t)}{t}=0$,
\item[(f4'')] $\lim_{t \to +\infty} \frac{f(t)}{t^{2^*-1}}=1$, where $2^*:=\frac{2N}{N-2}$,
\item[(f5'')] for some $C>0$ and $\max\{2^*-2, 2\} <p< 2^*$, i.e. satisfying
\begin{equation}\label{eq_cond_p_local}
p \in \parag{&(4, 6)& \quad N =3, \\ &\Big(2, \frac{2N}{N-2}\Big)& \quad N\geq 4, }
\end{equation}
(see also Remark \ref{rem_C_grande}), it results that
$$f(t) \geq t^{2^*-1} + C t^{p-1} \quad \textnormal{for $t \geq 0$}.$$
\end{itemize}

\begin{theorem}\label{teo_main_locale}
Suppose $N\geq 3$ and that \textnormal{(V1)-(V2)}, \textnormal{(f1'')--(f5'')} hold. 		
Let $K$ be defined by \eqref{eq_def_K}.
Then, up to a discretized subsequence, for small $\eps>0$ the equation 
\begin{equation}\label{eq_principale_locale}
-\eps^{2} \Delta u + V(x) u = f(u), \quad x \in \R^N
\end{equation}
has at least $\cupl(K)+1$ positive solutions, which belong to $C^{1, \sigma}(\R^N) \cap L^{\infty}(\R^N)$ for any $\sigma \in (0,1)$. Moreover, each of these sequences $u_{\eps}\in H^1(\R^N)$ concentrate in $K$ as $\eps \to 0$. Namely, for each small $\eps >0$ there exists a maximum point $x_{\eps}\in \R^N$ of $u_{\eps}$ such that 
$$\lim_{\eps \to 0}d(x_{\eps}, K) =0.$$
In addition, $ u_{\eps}(\eps \cdot+x_{\eps})$ converges in $H^1(\R^N)$ and uniformly on compact sets to a least energy solution of
\begin{equation}\label{eq_limite_locale}
-\Delta U + m_0 U = f(U), \quad x \in \R^N
\end{equation}
and, for some positive $C', C''$ independent on $\eps$, we have the uniform exponential decay
 $$ u_{\eps}(x) \leq C' \textnormal{exp}\Big(-C''\Big|\frac{x-x_{\eps}}{\eps}\Big|\Big) \quad\textit{ for $x \in \R^N$}.$$
\end{theorem}

\claim Proof.
The arguments of the previous sections apply mutatis mutandis. 
Indeed, we define in the same way the set of ground states $\widehat{S}$, which turns to be nonempty, compact and uniformly bounded in $L^{\infty}(\R^N)$ (see also \cite[Theorem 1.1]{ZhZo1}, \cite[Section 2]{BZZ} and \cite[Proposition 2.1]{ZCZ} for details). Then the truncation machinery can be set in motion, and one can prove $\widehat{S}_k \subset \widehat{S}$ in the same way as in Proposition \ref{prop_uguagl_insiem} (see \cite{BL, JT1, ZhZo1}). Existence, multiplicity and decay of solutions are given by \cite[Theorem 1.1 and Remark 1.3]{CJT}; the regularity is instead a consequence of standard elliptic estimates \cite[Appendix B]{Str}.
\QED

\bigskip

\medskip

\noindent {\bf Acknowledgments.} The author is supported by MIUR-PRIN project ``Qualitative and quantitative aspects of nonlinear PDEs'' (2017JPCAPN\_005), and partially supported by GNAMPA-INdAM. 
Moreover, the author would like to thank Prof. S. Cingolani for some fruitful comments.


\bigskip

\begin{thebibliography} {l}

 \bibitem{ADS} C.O. Alves, J.M. do Ó, M.A.S. Souto,
\emph{Local mountain-pass for a class of elliptic problems in $\R^N$ involving critical growth},
Nonlinear Anal. \textbf{46} (2001), no. 4, 495--510.

\bibitem{AlMi1} C.O. Alves, O.H. Miyagaki,
\emph{Existence and concentration of solution for a class of fractional elliptic equation in $\R^N$ via penalization method},
Calc. Var. Partial Differential Equations \textbf{55} (2016), no. 3, article ID 47, pp. 19. 
 
\bibitem{AlSo} C.O. Alves, M.A.S. Souto,
\emph{On existence and concentration behavior of ground state solutions for a class of problems with critical growth},
Comm. Pure Appl. Math. \textbf{1} (2002), no. 3, 417--431.

\bibitem{ABC} A. Ambrosetti, M. Badiale, S. Cingolani, 
\emph{Semiclassical states of nonlinear Schr\"odinger equations}, 
Arch. Ration. Mech. Anal. \textbf{140} (1997), no. 3, 285--300. 
	
\bibitem{AMS} A. Ambrosetti, A. Malchiodi, S. Secchi,
\emph{Multiplicity results for some nonlinear Schr\"odinger equations with potentials}, 
Arch. Ration. Mech. Anal. \textbf{159} (2001), no. 3, 253--271. 	
	
\bibitem{Amb4} V. Ambrosio,
\emph{Existence and concentration results for some fractional Schr\"odinger equations in $\R^N$ with magnetic fields},
Comm. Partial Differential Equations \textbf{44} (2019), no. 8, 637--680. 

\bibitem{Bar} T. Bartsch,
``Topological methods for variational problems with symmetries'',
Lectures Notes in Math. 1560, Springer-Verlag Berlin Heidelberg (1993).

\bibitem{BGM1} V. Benci, M. Ghimenti, A.M. Micheletti,
\emph{The nonlinear Schroedinger equation: solitons dynamics},
J. Differential Equations \textbf{249} (2010), no. 12, 3312--3341.
 
\bibitem{BGM2} V. Benci, M. Ghimenti, A.M. Micheletti,
\emph{On the dynamics of solitons in the nonlinear Schr\"odinger equation},
Arch. Ration. Mech. Anal. \textbf{205} (2012), no. 2, 467--492.
	 
\bibitem{BL} H. Berestycki, P.-L. Lions,
\emph{Nonlinear scalar field equations I: existence of a ground state},
Arch. Ration. Mech. Anal. \textbf{82} (1983), no. 4, 313--345.
	
\bibitem{BrJe} J.C. Bronski, R.L. Jerrard,
\emph{Soliton dynamics in a potential},
Math. Res. Lett. \textbf{7} (2002), no. 2-3, 329--342.
 
\bibitem{BJ3} J. Byeon, L. Jeanjean,
\emph{Standing waves for nonlinear Schr\"odinger equations with a general nonlinearity},
Arch. Ration. Mech. Anal. \textbf{185} (2007), no. 2, 185--200, and 
\emph{Erratum},
Arch. Ration. Mech. Anal. \textbf{190} (2008), no. 13, 549--551.

\bibitem{BKS} J. Byeon, O. Kwon, J. Seok,
\emph{Nonlinear scalar field equations involving the fractional Laplacian},
Nonlinearity \textbf{30} (2017), no. 4, 1659--1681.

\bibitem{BT1} J. Byeon, K. Tanaka,
\emph{Semi-classical standing waves for nonlinear Schr\"odinger equations at structurally stable critical points of the potential},
J. Eur. Math. Soc. (JEMS) \textbf{15} (2013), no. 5, 1859--1899. 

\bibitem{BZZ} J. Byeon, J. Zhang, W. Zou,
\emph{Singularly perturbed nonlinear Dirichlet problems involving critical growth},
Calc. Var. Partial Differential Equations \textbf{47} (2013), no. 1-2, 65--85.

\bibitem{Cha} K.-C. Chang,
``Infinite dimensional Morse theory and multiple solution problems'',
Progr. Nonlinear Differential Equations Appl. \textbf{6}, Birkh\"auser, Boston (1993).

\bibitem{Che} G. Chen,
\emph{Multiple semiclassical standing waves for fractional nonlinear Schr\"odinger equations},
Nonlinearity \textbf{28} (2015), no. 4, 927--949. 

\bibitem{CG} S. Cingolani, M. Gallo,
\emph{On the fractional NLS equation and the effects of the potential well’s topology},
Adv. Nonlinear Stud. (2021), DOI: \href{https://doi.org/10.1515/ans-2020-2114}{10.1515/ans-2020-2114}, pp. 40.

\bibitem{CJT} S. Cingolani, L. Jeanjean, K. Tanaka,
\emph{Multiplicity of positive solutions of nonlinear Schr\"odinger equations concentrating at a potential well},
Calc. Var. Partial Differential Equations \textbf{53} (2015), no. 1-2, 413--439.

\bibitem{CL} S. Cingolani, M. Lazzo, 
\emph{Multiple semiclassical standing waves for a class of nonlinear Schr\"odinger equations}, 
Topol. Methods Nonlinear Anal. \textbf{10} (1997), no. 1, 1--13. 

\bibitem{CLOT} O. Cornea, G. Lupton, J. Oprea, D. Tanr\'e,
``Lusternik-Schnirelmann Category'',
Math. Surveys Monogr. \textbf{103}, AMS (2003).
	
\bibitem{Coz} M. Cozzi,
\emph{Regularity results and Harnack inequalities for minimizers and solutions of nonlocal problems: a unified approach via fractional De Giorgi classes},
J. Funct. Anal. \textbf{272} (2017), no. 11, 4762--4837.

\bibitem{DPW} J. D\'avila, M. del Pino, J. Wei,
\emph{Concentrating standing waves for the fractional nonlinear Schr\"odinger equation},
J. Differential Equations \textbf{256} (2014), no. 2, 858--892.

\bibitem{DF} M. del Pino, P.L. Felmer, 
\emph{Local mountain passes for semilinear elliptic problems in unbounded domains},
Calc. Var. Partial Differential Equations \textbf{4} (1996), no. 2, 121--137. 
	
\bibitem{DGV} E. Di Nezza, G. Palatucci, E. Valdinoci,
\emph{Hitchhiker's guide to the fractional Sobolev spaces},
Bull. Sci. Math. \textbf{136} (2012), no. 5, 521--573.
 
\bibitem{DMV} S. Dipierro, M. Medina, E. Valdinoci,
``Fractional elliptic problems with critical growth in the whole of $\R^n$'',
Lecture Notes, Edizioni della Normale (2017), pp. 162.

\bibitem{FMV} M.M. Fall, F. Mahmoudi, E. Valdinoci,
\emph{Ground states and concentration phenomena for the fractional Schr\"odinger equation},
Nonlinearity \textbf{28} (2015), no. 6, 1937--1961.

\bibitem{FS} G. M. Figueiredo, G. Siciliano,
\emph{A multiplicity result via Ljusternick-Schnirelmann category and Morse theory for a fractional Schr\"odinger equation in $\R^N$},
NoDEA Nonlinear Differential Equations Appl. \textbf{23} 
(2016), no. 2, Article ID 12, pp. 22.

\bibitem{FW} A. Floer, A. Weinstein, 
\emph{Nonspreading wave packets for the cubic Schr\"odinger equation with a bounded potential},
J. Funct. Anal. \textbf{69} (1986), no. 3, 397--408. 

\bibitem{FLRW} G. Fournier, D. Lupo, M. Ramos, M. Willem,
\emph{Limit relative category and critical point theory},
in ``Dynamics reported. Expositions in dynamical systems \textbf{3} (New Series)'' (eds. C.K.R.T. Jones, U. Kirchgraber, H.O. Walther), Springer-Verlag, Berlin (1994), 1--24.

\bibitem{FGJS} J. Fr\"ohlich, S. Gustafson, B.L.G. Jonsson, I.M. Sigal,
\emph{Solitary wave dynamics in an external potential},
Comm. Math. Phys. \textbf{250} (2004), no. 3, 613--642.

\bibitem{Hat} A. Hatcher,
``Algebraic Topology'',
Cambridge Univ. Press (2002).

\bibitem{He2} Y. He,
\emph{Singularly perturbed fractional Schr\"odinger equations with critical growth},
Adv. Nonlinear Stud. \textbf{18} (2018), no. 3, pp. 25.

\bibitem{HeZo} Y. He, W. Zou,
\emph{Existence and concentration result for the fractional Schr\"odinger equations with critical nonlinearities},
Calc. Var. Partial Differential Equations \textbf{55} (2016), no. 4, 55--91.

\bibitem{JT1} L. Jeanjean, K. Tanaka,
\emph{A remark on least energy solutions in $\R^N$},
Proc. Amer. Math. Soc. \textbf{131} (2003), no. 8, 2399--2408.

\bibitem{JLZ} H. Jin, W. Liu, J. Zhang,
\emph{Singularly perturbed fractional Schr\"odinger equation involving a general critical nonlinearity},
Adv. Nonlinear Stud. \textbf{18} (2018), no. 3, pp. 13.

\bibitem{JFG} B.L.G. Jonsson, J. Fr\"ohlich, S. Gustafson,
\emph{Long time motion of NLS solitary waves in a confining potential},
Ann. Henri Poincar\'e \textbf{7} (2006), no. 4, 621--660.

\bibitem{Las} N. Laskin, 
\emph{Fractional quantum mechanics and L\'evy path integrals},
Phys. Lett. A \textbf{268} (2000), no. 4--6, 298--305.

\bibitem{Las2} N. Laskin, 
\emph{Fractional quantum mechanics},
Phys. Rev. E \textbf{62} (2000), no. 3, 3135--3145.

\bibitem{Las3} N. Laskin, 
\emph{Fractional Schr\"odinger equation},
Phys. Rev. E \textbf{66} (2002), no. 5, article ID 056108, pp. 7.

\bibitem{LTZW} Q. Li, K. Teng, J. Zhang, W. Wang,
\emph{Concentration behavior of solutions for fractional Schr\"odinger equations involving critical exponent},
J. Math. Phys. \textbf{61} (2020), no. 7, article ID 071513, pp. 5.

\bibitem{Oh2} Y.-G. Oh, 
\emph{Existence of semiclassical bound states of nonlinear Schr\"odinger equations with potentials of the class $(V)_a$},
Comm. Math. Phys. \textbf{13} (1988), no. 12, 1499--1519. 
	
\bibitem{Rab} P.H. Rabinowitz, 
\emph{On a class of nonlinear Schr\"odinger equations},
Z. Angew. Math. Phys. \textbf{43} (1992), no. 2, 270--291. 
	
\bibitem{SS} S. Secchi, M. Squassina,
\emph{Soliton dynamics for fractional Schr\"odinger equations},
Appl. Anal. \textbf{93} (2014), no. 8, 102--1729.

\bibitem{Seo} S. Seok,
\emph{Spike-layer solutions to nonlinear fractional Schr\"odinger equations with almost optimal nonlinearities},
Electron. J. Differential Equations, \textbf{2015} (2015), no. 196, 1--19.

\bibitem{SZ2} X. Shang, J. Zhang,
\emph{Existence and concentration of positive solutions for fractional nonlinear Schr\"odinger equation with critical growth},
J. Math. Phys. \textbf{58} (2017), no. 8, article ID 081502, pp. 18.

\bibitem{SZY} X. Shang, J. Zhang, Y. Yang,
\emph{On fractional Schr\"odinger equation in $\R^N$ with critical growth},
J. Math. Phys. \textbf{54} (2013), no. 12, article ID 121502, pp. 20.

\bibitem{Str} M. Struwe,
``Variational methods: applications to nonlinear partial differential equations and Hamiltonian systems'',
A series of modern surveys in mathematics \textbf{34}, Springer-Verlag Berlin Heidelberg (2008).

\bibitem{WLXF} J. Wang, D. Lu, J. Xu, F. Zhang,
\emph{Multiple positive solutions for semilinear Schr\"odinger equations with critical growth in $\R^N$},
J. Math. Phys. \textbf{56} (2015), no. 4, article ID 041502, pp. 27.

\bibitem{Wan} X. Wang,
\emph{On concentration of positive bound states of nonlinear Schr\"odinger equations},
Commun. Math. Phys. \textbf{153} (1993), no. 2, 229--244. 

\bibitem{ZCZ} J. Zhang, Z. Chen, W. Zou,
\emph{Standing waves for nonlinear Schr\"odinger equations involving critical growth},
J. London Math. Soc. \textbf{90} (2014), no. 2, 827--844.

\bibitem{ZMS} J. Zhang, J.M. do Ó, M. Squassina,
\emph{Fractional Schr\"odinger-Poisson systems with a general subcritical or critical nonlinearity},
Adv. Nonlinear Stud. \textbf{15} (2016), no. 1, 15--30.
 
\bibitem{ZhZo1} J. Zhang, W. Zou, 
\emph{A Berestycki-Lions theorem revisited},
Commun. Contemp. Math. \textbf{14} (2012), no. 5, article ID 1250033, pp. 14.
 
 \bibitem{ZhZo2} J. Zhang, W. Zou, 
\emph{Solutions concentrating around the saddle points of the potential for critical Schr\"odinger equations},
Calc. Var. Partial Differential Equations \textbf{54} (2015), no. 4, 4119--4142.

 \end{thebibliography}
\end{document}